\documentclass[12pt, a4paper]{article}
\usepackage{amsmath,amssymb,amsbsy,amsfonts,amsthm,latexsym,
                          amsopn,amstext,amsxtra,euscript,amscd}

\begin{document}

\newtheorem{theorem}{Theorem}
\newtheorem{lemma}[theorem]{Lemma}
\newtheorem{claim}[theorem]{Claim}
\newtheorem{cor}[theorem]{Corollary}
\newtheorem{prop}[theorem]{Proposition}
\newtheorem{definition}{Definition}
\newtheorem{question}[theorem]{Open Question}
\newtheorem{remark}{Remark}

\title{A note on the least totient of a residue class}

\author{{M. Z. Garaev}
\\
\normalsize{Instituto de Matem{\'a}ticas}
\\
\normalsize{Universidad Nacional Aut\'onoma de M\'exico}
\\
\normalsize{Campus Morelia, Apartado Postal 61-3 (Xangari)}
\\
\normalsize{C.P. 58089, Morelia, Michoac{\'a}n, M{\'e}xico} \\
\normalsize{\tt garaev@matmor.unam.mx}\\
}

\date{}

\pagenumbering{arabic}

\maketitle

\begin{abstract}
Let $q$ be a large prime number, $a$ be any integer,  $\varepsilon$
be a fixed small positive quantity. Friedlander and
Shparlinksi~\cite{FSh} have shown that there exists a positive
integer $n\ll q^{5/2+\varepsilon}$ such that $\phi(n)$ falls into
the residue class $a \pmod q.$ Here, $\phi(n)$ denotes Euler's
function. In the present paper we improve this bound to $n\ll
q^{2+\varepsilon}.$
\end{abstract}


\section {Introduction}

Let $q$ denote a large prime number, $a$ be any integer. Let
$N(q,a)$ denote the smallest positive integer $n$ for which
$\phi(n)\equiv a\pmod q.$ The number $N(q,a)$ exists. Indeed, if
$a+1\equiv 0\pmod q$ then one can take $n=q.$ Otherwise, one can
take $n$ to be a prime from the arithmetical progression $a+1 \pmod
q.$

The problem of upper bound estimates for $N(q,a)$ has been a subject
of study of the work of Friedlander and Shparlinski~\cite{FSh}. In
the present paper we obtain a new upper bound for $N(q,a).$

Throughout, we use the notation $A\lesssim B$ or $B\gtrsim A$ to
indicate that $A\ll Bq^{\varepsilon}$ for any fixed $\varepsilon>0,$
where the implied constant may depend on $\varepsilon.$

\begin{theorem}
\label{thm:main1} For any prime $q$ and integer $a,$ we have $
N(q,a)\lesssim q^{2}.$
\end{theorem}

Theorem~\ref{thm:main1} improves the bound $N(q,a)\lesssim q^{5/2}$
of~\cite{FSh}.

\bigskip

In the opposite direction, the recent result of Friedlander and
Luca~\cite{FL} implies that there exists a sequence of arithmetical
progressions $a_k\pmod {m_k}$ with $m_k\to\infty$ as $k\to\infty$
such that $N(m_k, a_k)$ exists and
$$
\frac{\log N(m_k, a_k)}{\log m_k}\to\infty  \quad {\rm as} \quad
k\to \infty.
$$

\bigskip

Following~\cite{FSh}, one can look for a solution of the congruence
\begin{equation}
\label{eqn:phi=a} \phi(n)\equiv a \pmod q
\end{equation}
among numbers of the form $n=p_1p_2p_3$ with primes $p_1,p_2,p_3.$
Here one can take $p_1,p_2,p_3$ to be primes that run certain
disjoint intervals $I_1,I_2,I_3\subset [2,q).$ This
converts~\eqref{eqn:phi=a} to the congruence
$$
(p_1-1)(p_2-1)(p_3-1)\equiv a \pmod q,\quad (p_1,p_2,p_3)\in
I_1\times I_2\times I_3.
$$
If we define $I_1$ to be the set of primes of the interval $[1,
q^{1/2+0.1\varepsilon}],$ then, using Karatsuba's estimate for
character sums with shifted primes $p_1-1$ and his method of solving
multiplicative ternary problems, one can derive that for
$\gcd(a,q)=1$ the number of solutions of this congruence is
asymptotically equal to
$$
\frac{|I_1||I_2||I_3|}{q-1}+\frac{\theta}{q-1}|I_1|q^{-\delta}q\sqrt{|I_2||I_3|},\qquad
|\theta|<1.
$$
From this one obtains the upper bound $N(q,a)\lesssim q^{5/2}.$ In
the present paper, we aggregate to this consideration one
consequence of Huxley's refinement of the Hal\'asz-Montgomery method
for large values of Dirichlet polynomials. This allows to get the
improved upper bound for $N(q,a).$

Our present application of the theory of large value estimates can
be compared with Lemma 4 of Friedlander and Iwaniec~\cite{FrIw}.

\section{Character sums and large value estimates}

Let $\chi$ be a nonprincipal character modulo a prime $q,$ \, $k$ be
any integer with $\gcd(k,q)=1,$ \, $p$ be a prime variable,
\,$\varepsilon$ be a small positive quantity. When
$q^{1/2+\varepsilon}<L<q,$ from Karatsuba's estimate it follows that
$$
\sum_{L/2<p\le L}\chi(p+k) \ll L^{1-\delta},\qquad
\delta=\delta(\varepsilon)>0.
$$
To prove our theorem, we will combine this estimate with Huxley's
refinement of the Hal\'asz-Montgomery method for large value
estimates. A sufficient for our purposes form of it is as follows.

Let $a_n$ be numbers with $|a_n|\lesssim 1,$ let $0< V\le N$ and let
$R$ be the number of characters $\chi\pmod q$ for which
$$
\Bigl|\sum_{n=N+1}^{2N}a_n\chi(n)\Bigr|\ge V.
$$
Then Huxley's refinement implies that
$$
R\lesssim \frac{N^2}{V^2}+\frac{qN^4}{V^6},
$$
see Mongomery~\cite{Mont}, Huxley~\cite{Hux}, Huxley and
Jutila~\cite{HJ}, Jutila~\cite{J}. This estimate is nontrivial when
$V>N^{3/4}$ and $N<q.$ In the case $N\ge q$ one has $RV^2\lesssim
N^2;$ in the case $V\le N^{3/4}$ and $N<q$ one has $RV^6\le
N^3RV^2\lesssim qN^4.$

The above sum can be replaced  with its $\ell$-th moment, where
$\ell$ is a fixed positive integer, and then we have the bound
$$
R\lesssim \frac{N^{2\ell}}{V^{2\ell}}+\frac{qN^{4\ell}}{V^{6\ell}}.
$$

More generally, the results on large values of Dirichlet polynomials
deal with upper bounds for the number of pairs
$(\sigma_r+it_r,\chi_r),$ with $\sigma_r\ge 0$ and certain
conditions on $t_r$ and characters $\chi_r$ (not necessarily
distinct), for which
$$
\Bigl|\sum_{n=N+1}^{2N}a_n\chi_r(n)n^{-\sigma_r-it_r}\Bigr|\ge V.
$$
Such a general consideration is important in applications to zero
density problems for $\zeta(s)$ and $L(s,\chi).$ For further key
references, see Bourgain~\cite{B1}, Harman~\cite{Har},
Heath-Brown~\cite{HB}, Ivic~\cite{Iv}.

\section{Proof of Theorem~\ref{thm:main1}}

We can assume that $a$ is relatively prime to $q,$ since otherwise
the statement is trivial in view of $\phi(q^2)\equiv 0\pmod q.$

Let $0<\varepsilon<0.1$ be fixed. Let $k=[1/\varepsilon].$ Put
$N=q^{1/(4k-1)},\, N_1=q^{1/2+0.1\varepsilon}.$

Let $I_1$ denote the set of primes $p_1\in (N_1/2, N_1].$ For
$j=2,3,\ldots, 6k+1$ let $I_j$ denote the set of primes of the
interval $(2^{-j}N, 2^{-j+1}N].$ Then, for sufficiently large $q,$
$$
\frac{N_1}{\log N_1}\ll |I_1|\ll \frac{N_1}{\log N_1}, \quad
\frac{N}{\log q}\ll |I_{j}|\ll \frac{N}{\log q},\quad 2\le j\le
6k+1.
$$
Consider the congruence
\begin{equation}
\label{eqn:6k+1}
(p_1-1)(p_2-1)\ldots (p_{6k+1}-1)\equiv a\pmod q,
\quad p_j\in I_j, \quad 1\le j\le 6k+1.
\end{equation}
Note that the left hand side is equal to $\phi(p_1p_2\ldots
p_{6k+1})$ and
$$
p_1p_2\ldots p_{6k+1}\ll q^{6k/(4k-1)+1/2+0.1\varepsilon}\ll
q^{2+0.5\varepsilon}.
$$
Hence, since $N(q,a)$ exists, it suffices to prove that
congruence~\eqref{eqn:6k+1} has a solution for any sufficiently
large prime $q.$

Assume the contrary. We express the number of solutions of
congruence~\eqref{eqn:6k+1} (which is equal to zero by the
assumption) via character sums. Separating the contribution of the
principal character, we deduce
\begin{eqnarray*}
&& |I_1||I_2|\ldots
|I_{6k+1}| \le\\
&&\quad \sum_{\chi\not=\chi_0}\Bigl|\sum_{p_1\in
I_1}\chi(p_1-1)\Bigr|\Bigl|\sum_{p_2\in I_2}\chi(p_2-1)\Bigr|\ldots
\Bigl|\sum_{p_{6k+1}\in I_{6k+1}}\chi(p_{6k+1}-1)\Bigr|.
\end{eqnarray*}
The left hand side is $\gtrsim N^{6k}N_1.$ Hence, for some $2\le
j\le 6k+1$ and $I=I_j,$ we have
$$
N^{6k}N_1\lesssim \sum_{\chi\not=\chi_0}\Bigl|\sum_{p\in
I}\chi(p-1)\Bigr|^{6k}\Bigl|\sum_{p_1\in I_1}\chi(p_1-1)\Bigr|.
$$
Decomposing into level sets, for some positive numbers $V$ and $V_1$
we get that
\begin{equation}
\label{eqn:ErrorTh1} N^{6k}N_1\lesssim RV^{6k}V_1,
\end{equation}
where $R$ is the number of non-principal characters $\chi\pmod q$
for which
$$
V\le \Bigl|\sum_{p\in I}\chi(p-1)\Bigr|\le 2V,\quad V_1\le
\Bigl|\sum_{p_1\in I_1}\chi(p_1-1)\Bigr|\le 2V_1.
$$
By Karatsuba's estimate,
$$
V_1\ll N_1^{1-\delta},\quad \delta=\delta(\varepsilon)>0.
$$
From the large values estimate,
\begin{equation}
\label{eqn:LVE} R\lesssim
\frac{N^{6k}}{V^{6k}}+\frac{qN^{12k}}{V^{18k}}.
\end{equation}
Incorporating these estimates in~\eqref{eqn:ErrorTh1}, we get that
$$
N^{6k}N_1\lesssim
\Bigl(N^{6k}+\frac{qN^{12k}}{V^{12k}}\Bigr)N_1^{1-\delta}.
$$
Comparing the orders of the implied expressions, we obtain
$$
N^{6k}\lesssim \frac{qN^{12k}}{V^{12k}}N_1^{-\delta}.
$$
Therefore, from~\eqref{eqn:LVE} we get that
\begin{equation}
\label{eqn:18k}
RV^{18k}\lesssim
\Bigl(N^{6k}+\frac{qN^{12k}}{V^{12k}}\Bigr)V^{12k}\lesssim qN^{12k}.
\end{equation}
Since $N^{4k-1}=q,$ the congruence
$$
(x_1-1)\ldots (x_{4k-1}-1)\equiv (y_1-1)\ldots (y_{4k-1}-1)\pmod
q,\quad x_j, y_j\in I,
$$
implies the equality
$$
(x_1-1)\ldots (x_{4k-1}-1)=(y_1-1)\ldots (y_{4k-1}-1), \quad x_j,
y_j\in I.
$$
This equality has $\lesssim N^{4k-1}$ solutions. Hence,
\begin{equation}
\label{eqn:8k-2} RV^{8k-2}\le \sum_{\chi}\Bigl|\sum_{p\in
I}\chi(p-1)\Bigr|^{8k-2}\lesssim qN^{4k-1}.
\end{equation}
Let us estimate $RV_1^4.$ Recall that $N_1^2>q$ and observe that the
number of solutions of the congruence
$$
(x_1-1)(x_2-1)\equiv (x_3-1)(x_4-1)\pmod q,\quad x_j\in I_1,
$$
is not greater than twice the number of solutions of the equality
$$
(x_1-1)(x_2-1)=(x_3-1)(x_4-1)+tq, \qquad x_j\in I_1, \quad 0\le t\le
N_1^2/q.
$$
The right hand side of this equality does not vanish, so for each
triple $x_3, x_4, t$ we have $\lesssim 1$ choices for $x_1, x_2.$
Thus, the above congruence has $\lesssim N_1^4/q$ solutions. Hence,
\begin{equation}
\label{eqn:RV14} RV_1^4\le \sum_{\chi}\Bigl|\sum_{p_1\in
I_1}\chi(p_1-1)\Bigr|^4\lesssim N_1^4.
\end{equation}
Rewrite~\eqref{eqn:ErrorTh1} in the form
$$
N^{24k}N_1^{4}\, \lesssim
\,(RV^{18k})^{3/(5k+1)}(RV^{8k-2})^{15k/(5k+1)}RV_1^4.
$$
Taking into account the estimates~\eqref{eqn:18k}--\eqref{eqn:RV14},
we get that
$$
N\lesssim q^{(5k+1)/(20k^2+k)}.
$$
This contradicts to $N=q^{1/(4k-1)}.$

\bigskip

\section{Remark}

The reader may note, that in our treating of
congruence~\eqref{eqn:6k+1} the only essential property of the sets
$I_j$ \, $(j\ge 2)$ that we use is their density in the
corresponding intervals, and the set $I_1$ is used to apply the
nontrivial character sum estimate of Karatsuba. Thus, our argument
can be applied to deal with a class of other congruences.

\bigskip

{\bf Acknowledgement.} The author is grateful to S. V. Konyagin for
reading the paper.

\end{document}